\documentclass[preprint,12pt]{elsarticle}
\usepackage{etoolbox}
\makeatletter
\patchcmd{\ps@pprintTitle}{\footnotesize\itshape
       Preprint submitted to \ifx\@journal\@empty Elsevier
       \else\@journal\fi\hfill\today}{\relax}{}{}
\makeatother
\usepackage{graphicx}
\usepackage[section]{placeins}
 \usepackage{amssymb}
 \usepackage{lineno}
   \usepackage{amsthm}
   \usepackage{amsmath}
   \numberwithin{equation}{section}
     \begin{document}
\begin{frontmatter}
\title{New Two Step Laplace Adam-Bashforth Method for Integer an Non integer Order Partial Differential Equations}
\author{Rodrigue Gnitchogna*}
\author{Abdon Atangana**}
\address{*Department of Mathematics,\\ Faculty of Natural and Agricultural Science\\ University of the Free State Bloemfontein.}
\address{*Department of Mathematics,\\  University of Namibia.}
\address{**Institute of Groundwater Studies,\\ University of the Free State Bloemfontein South Africa.}

\begin{abstract}
This paper presents a novel method that allows to generalise the use of the Adam-Bashforth to Partial Differential Equations with local and non local operator.  The Method derives a two step Adam-Bashforth numerical scheme in Laplace space and the solution is taken back into the real space via inverse Laplace transform. The method yields a powerful numerical algorithm for fractional order derivative where the usually very difficult to manage summation in the numerical scheme disappears.  Error Analysis of the method is also presented. Applications of the method and numerical simulations are presented on a wave-equation like, and on a fractional order diffusion equation.
\end{abstract}
\end{frontmatter}
\section{Introduction}
Adam-Bashforth method has been recognised as a powerful numerical tool to solve Partial Differential equations ( P.D.E) [1]. It is a numerical scheme that is used  in many field of applied science. epidemiology, engineering in dynamical systems, in chaotic problems [2]-[6]. The numerical scheme is good for both differential equations with classical derivatives, and differential equations with non integer order derivatives. However The method is not fully applied to P.D.E with local  and non local operator as it was designed only for Ordinary Differential Equations (O.D.E) [1]. Various other methods are used instead for P.D.E with integer order differentiation, and those with real order derivatives[7]-[15].  Nonetheless, Due to the accuracy and efficiency of  Adam-Bashforth techniques, there is a need to extend the methods to P.D.E.  By eliminating one variable and transforming a P.D.E to an O.D.E via Laplace transform, the newly obtained O.D.E can be analysed in Laplace space, and a further application of the inverse transform will return the solution to the real space. This paper derives a new numerical scheme that combines the Laplace transform and Adam-Bashforth method to handle P.D.E with integer order and non integer order derivatives.

\section{ Numerical Method for P.D.E with Integer order}
Consider the following general P.D.E
\begin{equation}
\frac{\partial u(x,t)}{\partial t}=L u(x,t)+N u(x,t)
\end{equation}
Where $L$ is a linear operator and $N$ a non linear operator.
We start by applying Laplace transform on both sides of the equation, with respect to the variable $x$ to obtain

\[ \mathcal{L}\bigg(\frac{\partial u(x,t)}{\partial t}\bigg) =\mathcal{L}\bigg(L u(x,t)+N u(x,t)\bigg) \Longrightarrow\]
\[  \frac{d}{dt}(u(p,t))=\mathcal{L}\bigg(L u(x,t)+N u(x,t)\bigg) \]
\begin{equation}
 \frac{d}{dt}(u(t))=F(u,t) \label{eq2}
\end{equation}
Where $ u(t)=u(p,t)$ and $F(u,t)=\mathcal{L}\bigg(L u(x,t)+N u(x,t)\bigg). $

Next we apply the fundamental theorem of calculus on equation (\ref{eq2}) to obtain
\[ u(t)=u(t_{0})+\int_{0}^{t}{F(u,\tau) d\tau} \] this is also
\[ u(t)=u_{0}+\int_{0}^{t}{F(u,\tau) d\tau} \]
When $ t=t_{n+1}$ we have 
\[ u_{n+1}= u(t_{n+1})=u_{0}+\int_{0}^{t_{n+1}}{F(u,\tau) d\tau}    \]
When $t=t_{n}$
\[  u_{n}=u(t_{n})=u_{0}+\int_{0}^{t_{n}}{F(u,\tau) d\tau}   \] it follows that 
\[ u_{n+1}-u_{n}= \int_{0}^{t_{n+1}}{F(u,\tau) d\tau} -\int_{0}^{t_{n}}{F(u,\tau) d\tau} \]
\[ u_{n+1}-u_{n}= \int_{t_{n}}^{t_{n+1}}{F(u,\tau) d\tau} \]
If we approximate $F(u,t)$ with the Lagrange polynomial 
\[P(t)(\approx F(u,t))=  \frac{t-t_{n-1}}{t_{n}-t_{n-1}} F(u,t_{n})+\frac{t-t_{n}}{t_{n-1}-t_{n}} F(u,t_{n-1})\]
\[P(t)=  \frac{t-t_{n-1}}{t_{n}-t_{n-1}} F_{n}+\frac{t-t_{n}}{t_{n-1}-t_{n}} F_{n-1}\]
we can therefore write
 \[u_{n+1}-u_{n}=\int_{t_{n}}^{t_{n+1}}{F(u,\tau) d\tau} \]
 \[u_{n+1}-u_{n}=\int_{t_{n}}^{t_{n+1}}{\bigg(\frac{t-t_{n-1}}{t_{n}-t_{n-1}} F_{n}+\frac{t-t_{n}}{t_{n-1}-t_{n}} F_{n-1}\bigg) dt}\]
 \[u_{n+1}-u_{n}=\frac{F_{n}}{t_{n}-t_{n-1}}\int_{t_{n}}^{t_{n+1}}{\big({t-t_{n-1}} \big) dt}+\frac{F_{n-1}}{{t_{n-1}-t_{n}}}\int_{t_{n}}^{t_{n+1}}{\big({t_{n-1}-t_{n}}\big) dt}\]
 \[u_{n+1}-u_{n}=\frac{F_{n}}{t_{n}-t_{n-1}}\bigg[\frac{1}{2} t^{2}-t t_{n-1}\bigg]_{t_{n}}^{t_{n+1}}+\frac{F_{n-1}}{{t_{n-1}-t_{n}}}\bigg[\frac{1}{2}t^{2}-t t_{n}\bigg]_{t_{n}}^{t_{n+1}}\]
 By letting \[ h=t_{n}-t_{n-1}\] we have 
\begin{eqnarray*}u_{n+1}-u_{n}=\frac{F_{n}}{h}\bigg(\frac{1}{2} t_{n+1}^2-t_{n+1}t_{n-1}-\frac{1}{2}t_{n}^{2}+t_{n}t_{n-1} \bigg)-\frac{F_{n-1}}{h}\bigg(\frac{1}{2}t_{n+1}^{2}- t_{n}t_{n+1}-\frac{1}{2}t_{n}^2+t_{n}^{2}\bigg)\end{eqnarray*}
\[u_{n+1}-u_{n}=\frac{F_{n}}{h}\bigg(\frac{1}{2} (t_{n+1}-t_{n})(t_{n+1}+t_{n})-t_{n-1}(t_{n+1}-t_{n}) \bigg)-\frac{F_{n-1}}{h}\bigg(\frac{1}{2}(t_{n+1}- t_{n})(t_{n+1}+t_{n})-t_{n}(t_{n+1}-t_{n})\bigg)\]
\begin{eqnarray*}u_{n+1}-u_{n}=\frac{F_{n}}{h}\bigg(\frac{1}{2} h(t_{n+1}+t_{n})-ht_{n-1}\bigg)-\frac{F_{n-1}}{h}\bigg(\frac{1}{2}h(t_{n+1}+t_{n})-ht_{n}\bigg)\end{eqnarray*}
\begin{eqnarray*}u_{n+1}-u_{n}=F_{n}\bigg(\frac{1}{2} (t_{n+1}+t_{n})-t_{n-1}\bigg)-F_{n-1}\bigg(\frac{1}{2}(t_{n+1}+t_{n})-t_{n}\bigg)\end{eqnarray*}
\begin{eqnarray*}u_{n+1}-u_{n}=F_{n}\bigg(\frac{1}{2}\Big((n+1)h+{n}h\Big)-(n-1)h\bigg)-F_{n-1}\bigg(\frac{1}{2}\Big(({n+1})h+{n}h\Big)-{n}h\bigg)\end{eqnarray*}
\begin{eqnarray*}u_{n+1}-u_{n}=F_{n}\bigg(nh+\frac{1}{2}h-n h+h\bigg)-F_{n-1}\bigg(n h+\frac{1}{2}h-n h\bigg)\end{eqnarray*}
\begin{eqnarray}u_{n+1}=u_{n}+h\bigg(\frac{3}{2}F_{n}-\frac{1}{2}F_{n-1}\bigg)  \label{eq23}\end{eqnarray}
Applying the inverse Laplace transform to return into the real space, we have:
\[ \mathcal{L}^{-1}(u_{n+1})=\mathcal{L}^{-1}\bigg[u_{n}+h\bigg(\frac{3}{2}F_{n}-\frac{1}{2}F_{n-1}\bigg)\bigg] \]
\[ u(x,t)=\mathcal{L}^{-1}\bigg[u_{n}+h\bigg(\frac{3}{2}F_{n}-\frac{1}{2}F_{n-1}\bigg)\bigg] \]

\section{New Numerical Method for P.D.E with non integer order}
To illustrate the method we consider the general fractional P.D.E 
\begin{equation}
\frac{\partial^{\alpha} u(x,t)}{\partial t^{\alpha}}=L u(x,t)+N u(x,t) \label{eq3}
\end{equation}
Where $L$ is a linear operator and $N$ a non linear operator.
Applying Laplace transform on both sides of the equation (\ref{eq3}), we have
\[ \mathcal{L}\Big( \frac{\partial^{\alpha} u(x,t)}{\partial t^{\alpha}}  \Big)=\mathcal{L}\Big(L u(x,t)+N u(x,t)   \Big) \]
For the Caputo type fractional partial derivative this will be
\[ {}^{C}_{a}\!\textit{D}^{\alpha}_{t}u(p,t) = \mathcal{L}\Big(L u(x,t)+N u(x,t)   \Big)\]
\[ {}^{C}_{a}\!\textit{D}^{\alpha}_{t}u(p,t) =F(u,t)   \Big)\] this is
\begin{equation}
{}^{C}_{a}\!\textit{D}^{\alpha}_{t}u(t) =F(u,t)  \label{eq32}
\end{equation}
where $ u(t)=u(p,t) \qquad\textit{and} \qquad F(u,t)=\mathcal{L}\Big(L u(x,t)+N u(x,t) \Big)$.
\\The next step is to apply the Caputo fractional Integral operator on equation (\ref{eq32}). Doing so we obtain
\[u(t)-u(t_0)=\frac{1}{\Gamma(\alpha)} \int_{0}^{t} {(t-\tau)^{\alpha-1}F(u,\tau)d\tau} \]
When $t=t_{n+1}$
\begin{eqnarray*}u_{n+1}=u(t_{n+1})=u_{0}+\frac{1}{\Gamma(\alpha)} \int_{0}^{t_{n+1}}{(t_{n+1}-\tau)^{\alpha-1}F(u,\tau)d\tau} \end{eqnarray*}
When $t=t_{n}$
\begin{eqnarray*}u_{n}=u(t_{n})=u_{0}+\frac{1}{\Gamma(\alpha)} \int_{0}^{t_{n}} {(t_{n}-\tau)^{\alpha-1}F(u,\tau)d\tau} \end{eqnarray*}
\begin{eqnarray}u_{n+1}-u_{n}=\frac{1}{\Gamma(\alpha)} \Big[\int_{0}^{t_{n+1}}{(t_{n+1}-\tau)^{\alpha-1}F(u,\tau)d\tau-\int_{0}^{t_{n}} {(t_{n}-\tau)^{\alpha-1}F(u,\tau)d\tau}}\Big]\qquad \label{eq33} \end{eqnarray}
\[\int_{0}^{t_{n+1}}{(t_{n+1}-\tau)^{\alpha-1}F(u,\tau)d\tau}=\sum_{j=0}^{n}{\int_{t_{j}}^{t_{j+1}}(t_{n+1}-\tau)^{\alpha-1}F(u,\tau)d\tau} \]
We approximate $F(u,t)$ with the following Lagrange polynomial
\[P(t)(\approx F(u,t))=  \frac{t-t_{n-1}}{t_{n}-t_{n-1}} F(u,t_{n})+\frac{t-t_{n}}{t_{n-1}-t_{n}} F(u,t_{n-1})\]
\[P(t)=  \frac{t-t_{n-1}}{t_{n}-t_{n-1}} F_{n}+\frac{t-t_{n}}{t_{n-1}-t_{n}} F_{n-1}\]

The first fractional integral in equation (\ref{eq33}) can then be expressed as

\[\int_{0}^{t_{n+1}}{(t_{n+1}-\tau)^{\alpha-1}F(u,\tau)d\tau}=\sum_{j=0}^{n}{\int_{t_{j}}^{t_{j+1}}(t_{n+1}-t)^{\alpha-1}\Big(\frac{t-t_{n-1}}{t_{n}-t_{n-1}} F_{n}+\frac{t-t_{n}}{t_{n-1}-t_{n}} F_{n-1}\Big)dt} \]
\[=\sum_{j=0}^{n}\bigg[\frac{F_{n}}{t_{n}-t_{n-1}} \int_{t_{j}}^{t_{j+1}}(t_{n+1}-t)^{\alpha-1}({t-t_{n-1}})dt+\frac{F_{n-1}}{t_{n-1}-t_{n}}\int_{t_{j}}^{t_{j+1}}(t_{n+1}-t)^{\alpha-1}({t-t_{n}})dt\bigg] \]
\[=\sum_{j=0}^{n}\bigg[\frac{F_{n}}{h} \int_{t_{j}}^{t_{j+1}}(t_{n+1}-t)^{\alpha-1}({t-t_{n-1}})dt-\frac{F_{n-1}}{h}\int_{t_{j}}^{t_{j+1}}(t_{n+1}-t)^{\alpha-1}({t-t_{n}})dt\bigg] \]
We now implement the following change of variable. \\We let $y=t_{n+1}-t, dt=-dy, t=t_{n+1}-y$
\begin{eqnarray*}  \int_{t_{j}}^{t_{j+1}}(t_{n+1}-t)^{\alpha-1}({t-t_{n-1}})dt =\int_{t_{n+1}-t_{j}}^{t_{n+1}-t_{j+1}}y^{\alpha-1}(-y+{t_{n+1}-t_{n-1}})dy\end{eqnarray*}
\[ =\int_{t_{n+1}-t_{j}}^{t_{n+1}-t_{j+1}}(y^{\alpha}-2hy^{\alpha-1})dy\]
\[ =\frac{1}{\alpha+1}\Big[y^{\alpha+1}\Big]_{t_{n+1}-t_{j}}^{t_{n+1}-t_{j+1}}-\frac{2h}{\alpha}\Big[y^{\alpha}\Big]_{t_{n+1}-t_{j}}^{t_{n+1}-t_{j+1}}    \]
\[=\frac{1}{\alpha+1} \Big((t_{n+1}-t_{j+1})^{\alpha+1}-(t_{n+1}-t_{j})^{\alpha+1}\Big)-\frac{2h}{\alpha}\Big(  (t_{n+1}-t_{j+1})^{\alpha}-(t_{n+1}-t_{j})^{\alpha} \Big)\].
On the other hand
\[\int_{t_{j}}^{t_{j+1}}(t_{n+1}-t)^{\alpha-1}({t-t_{n}})dt=-\int_{t_{n+1}-t_{j}}^{t_{n+1}-t_{j+1}}y^{\alpha-1}({-y+t_{n+1}-t_{n}})dy\]
\[=\int_{t_{n+1}-t_{j}}^{t_{n+1}-t_{j+1}}(y^{\alpha}-hy^{\alpha-1})dy\]
\[ =\frac{1}{\alpha+1}\Big[y^{\alpha+1}\Big]_{t_{n+1}-t_{j}}^{t_{n+1}-t_{j+1}}-\frac{h}{\alpha}\Big[y^{\alpha}\Big]_{t_{n+1}-t_{j}}^{t_{n+1}-t_{j+1}}    \]
\[=\frac{1}{\alpha+1} \Big((t_{n+1}-t_{j+1})^{\alpha+1}-(t_{n+1}-t_{j})^{\alpha+1}\Big)-\frac{h}{\alpha}\Big(  (t_{n+1}-t_{j+1})^{\alpha}-(t_{n+1}-t_{j})^{\alpha} \Big)\]
It will follow that 
\begin{eqnarray*}
\int_{0}^{t_{n+1}}{(t_{n+1}-\tau)}^{\alpha-1}F(u,\tau)d\tau \end{eqnarray*}
\begin{eqnarray*}
=\frac{F_{n}}{h}\bigg\{ \sum_{j=0}^{n}\bigg[ \frac{1}{\alpha+1} \Big((t_{n+1}-t_{j+1})^{\alpha+1}-(t_{n+1}-t_{j})^{\alpha+1}\Big)\bigg]\\
-\frac{2h}{\alpha}\sum_{j=0}^{n}\bigg[\Big((t_{n+1}-t_{j+1})^{\alpha}-(t_{n+1}-t_{j})^{\alpha} \Big)\bigg]\bigg\}\end{eqnarray*}
\begin{eqnarray*}
-\frac{F_{n-1}}{h}\bigg\{ \sum_{j=0}^{n}\bigg[ \frac{1}{\alpha+1} \Big((t_{n+1}-t_{j+1})^{\alpha+1}-(t_{n+1}-t_{j})^{\alpha+1}\Big)\bigg]\\
-\frac{h}{\alpha}\sum_{j=0}^{n}\bigg[\Big((t_{n+1}-t_{j+1})^{\alpha}-(t_{n+1}-t_{j})^{\alpha} \Big)\bigg]\bigg\}\end{eqnarray*}
\begin{eqnarray*}
=\frac{F_{n}}{h}\bigg( \frac{1}{\alpha+1}(-(t_{n+1}-t_{0})^{\alpha+1})-\frac{2h}{\alpha}(-(t_{n+1}-t_{0})^{\alpha})\bigg)\end{eqnarray*}
\begin{eqnarray*}
-\frac{F_{n-1}}{h}\bigg(\frac{1}{\alpha+1} (-(t_{n+1}-t_{0})^{\alpha+1})-\frac{h}{\alpha}(-(t_{n+1}-t_{0})^{\alpha})\bigg)\end{eqnarray*}
\begin{eqnarray*}
=\frac{F_{n}}{h}\bigg( \frac{-(n+1)^{\alpha+1}h^{\alpha+1}}{\alpha+1}+\frac{2h ({n+1})^{\alpha}h^{\alpha}}{\alpha}\bigg)\\
-\frac{F_{n-1}}{h}\bigg(\frac{-(n+1)^{\alpha+1}h^{\alpha+1}}{\alpha+1} +\frac{h({n+1})^{\alpha}h^{\alpha}}{\alpha}\bigg).
\end{eqnarray*}
\begin{eqnarray*}
=h^{\alpha}\Big[\bigg( \frac{2 ({n+1})^{\alpha}}{\alpha}-\frac{(n+1)^{\alpha+1}}{\alpha+1}\bigg)F_{n}&\\
-\bigg(\frac{(n+1)^{\alpha}}{\alpha}-\frac{({n+1})^{\alpha+1}}{\alpha+1}\bigg)F_{n-1}\Big].
\end{eqnarray*}
The second fractional integral in equation (\ref{eq33}) will similarly be evaluated as 
\[\int_{0}^{t_{n}}\!\!\!\!\!\!{(t_{n}-\tau)^{\alpha-1}F(u,\tau)d\tau}=\sum_{j=0}^{n-1}{\int_{t_{j}}^{t_{j+1}}(t_{n}-t)^{\alpha-1}\Big(\frac{t-t_{n-1}}{t_{n}-t_{n-1}} F_{n}+\frac{t-t_{n}}{t_{n-1}-t_{n}} F_{n-1}\Big)dt} \]
\[\int_{0}^{t_{n}}\!\!\!\!\!\!{(t_{n}-\tau)^{\alpha-1}F(u,\tau)d\tau}=\frac{F_{n}}{h}\sum_{j=0}^{n-1}{\int_{t_{j}}^{t_{j+1}}\!\!\!\!\!\!\!\!(t_{n}-t)^{\alpha-1}({t-t_{n-1}})dt} -\frac{F_{n-1}}{h}\sum_{j=0}^{n-1}{\int_{t_{j}}^{t_{j+1}}\!\!\!\!\!\!\!\!(t_{n}-t)^{\alpha-1}\!({t-t_{n}})dt} \]
We also implement the following change of variable. \\We let $y=t_{n}-t, dt=-dy, t=t_{n}-y$
\[\int_{0}^{t_{n}}\!\!\!\!\!\!{(t_{n}-\tau)^{\alpha-1}F(u,\tau)d\tau}=\frac{F_{n}}{h}\sum_{j=0}^{n-1}{\int_{t_{n}-t_{j}}^{t_{n}-t_{j+1}}\!\!\!\!\!\!\!\!-(y)^{\alpha-1}(t_{n}-t_{n-1}-y)dy}-\frac{F_{n-1}}{h}\sum_{j=0}^{n-1}{\int_{t_{n}-t_{j}}^{t_{n}-t_{j+1}}y^{\alpha}dy} \]
\[\int_{0}^{t_{n}}\!\!\!\!\!\!{(t_{n}-\tau)^{\alpha-1}F(u,\tau)d\tau}=\frac{F_{n}}{h}\sum_{j=0}^{n-1}{\int_{t_{n}-t_{j}}^{t_{n}-t_{j+1}}\!\!\!\!\!\!\!\!(y^{\alpha}-hy^{\alpha-1})dy}-\frac{F_{n-1}}{h}\sum_{j=0}^{n-1}{\int_{t_{n}-t_{j}}^{t_{n}-t_{j+1}}y^{\alpha}dy} \]
\[ \int_{0}^{t_{n}}\!\!\!\!\!\!{(t_{n}-\tau)^{\alpha-1}F(u,\tau)d\tau}=\frac{F_{n}}{h}\sum_{j=0}^{n-1}{\Big[\frac{y^{\alpha+1}}{\alpha+1}-\frac{h}{\alpha}y^{\alpha}\Big]}_{t_{n}-t_{j}}^{t_{n}-t_{j+1}}-\frac{F_{n-1}}{h}\sum_{j=0}^{n-1}{ \Big[\frac{y^{\alpha+1}}{\alpha+1}\Big]}_{t_{n}-t_{j}}^{t_{n}-t_{j+1}}\]
\begin{eqnarray*} =\frac{F_{n}}{h}\sum_{j=0}^{n-1}{\Big(\frac{(t_{n}-t_{j+1})^{\alpha+1}}{\alpha+1}-\frac{h}{\alpha}(t_{n}-t_{j+1})^{\alpha}-\frac{(t_{n}-t_{j})^{\alpha+1}}{\alpha+1}+\frac{h}{\alpha}(t_{n}-t_{j})^{\alpha}\Big)}\\
-\frac{F_{n-1}}{h}\sum_{j=0}^{n-1}{ \Big(\frac{(t_{n}-t_{j+1})^{\alpha+1}}{\alpha+1}-\frac{(t_{n}-t_{j})^{\alpha+1}}{\alpha+1}\Big)}\end{eqnarray*}
\begin{eqnarray*} =\frac{F_{n}}{h}\Bigg\{\sum_{j=0}^{n-1}{\Big(\frac{(t_{n}-t_{j+1})^{\alpha+1}}{\alpha+1}-\frac{(t_{n}-t_{j})^{\alpha+1}}{\alpha+1}\Big)}-\frac{h}{\alpha}\sum_{j=0}^{n-1}{\Big((t_{n}-t_{j+1})^{\alpha}-(t_{n}-t_{j})^{\alpha}\Big)}\Bigg\}\\
-\frac{F_{n-1}}{h}\sum_{j=0}^{n-1}{ \Big(\frac{(t_{n}-t_{j+1})^{\alpha+1}}{\alpha+1}-\frac{(t_{n}-t_{j})^{\alpha+1}}{\alpha+1}\Big)}\end{eqnarray*}
\begin{eqnarray*} =\frac{F_{n}}{h}{\Big\{-\frac{(t_{n}-t_{0})^{\alpha+1}}{\alpha+1}}-\frac{h}{\alpha}{\Big(-(t_{n}-t_{0})^{\alpha}\Big)}\Big\}-\frac{F_{n-1}}{h(\alpha+1)}\big(-(t_{n}-t_{0})^{\alpha+1}\big)\end{eqnarray*}
\begin{eqnarray*} =\frac{F_{n}}{h}{\Big(-\frac{(t_{n}-t_{0})^{\alpha+1}}{\alpha+1}}+\frac{h}{\alpha}{(t_{n}-t_{0})^{\alpha}}\Big)+\frac{F_{n-1}}{h(\alpha+1)}(t_{n}-t_{0})^{\alpha+1}\end{eqnarray*}
\begin{eqnarray*} =\frac{F_{n}}{h}{\Big(\frac{-n^{\alpha+1} h^{\alpha+1}}{\alpha+1}}+\frac{ n^{\alpha} h^{\alpha+1}}{\alpha}\Big)+\frac{n^{\alpha+1} h^{\alpha+1}}{h(\alpha+1)}F_{n-1}\end{eqnarray*}
 Therefore we can then write 
\begin{eqnarray*}
\int_{0}^{t_{n}}\!\!\!\!\!\!{(t_{n}-\tau)^{\alpha-1}F(u,\tau)d\tau}=h^{\alpha}{\bigg(\Big(\frac{ n^{\alpha} }{\alpha}-\frac{n^{\alpha+1} }{\alpha+1}}\Big) F_{n} +\frac{n^{\alpha+1}}{\alpha+1}F_{n-1}\bigg)
\end{eqnarray*}
Rewriting equation (\ref{eq33}) substituting in the later results we can then have
\begin{eqnarray*} u_{n+1}-u_{n}=\frac{h^{\alpha}}{\Gamma(\alpha)} \Big[\bigg( \frac{2 ({n+1})^{\alpha}}{\alpha}-\frac{(n+1)^{\alpha+1}}{\alpha+1}\bigg)F_{n}&\\
-\bigg(\frac{(n+1)^{\alpha}}{\alpha} -\frac{({n+1})^{\alpha+1}}{\alpha+1}\bigg)F_{n-1}-{\bigg(\Big(\frac{ n^{\alpha} }{\alpha}-\frac{n^{\alpha+1} }{\alpha+1}}\Big) F_{n} +\frac{n^{\alpha+1}}{\alpha+1}F_{n-1}\bigg) \Big]\end{eqnarray*}

\begin{eqnarray} u_{n+1}-u_{n}=\frac{h^{\alpha}}{\Gamma(\alpha)} \Big[\bigg( \frac{2 ({n+1})^{\alpha}-n^{\alpha}}{\alpha}+\frac{n^{\alpha+1}-(n+1)^{\alpha+1}}{\alpha+1}\bigg)F_{n}\label{eqnum1}&\\ 
-\bigg(\frac{(n+1)^{\alpha}}{\alpha} +\frac{n^{\alpha+1}-(n+1)^{\alpha+1}}{\alpha+1}\bigg)F_{n-1} \Big] \nonumber   \end{eqnarray}

\newtheorem{remark}{Remark}[section]
\begin{remark}
 For $\alpha =1$ we recover the classic Adam-Mashforth Numerical scheme. In fact 
 \[  u_{n+1}-u_{n}=\\\frac{h}{\Gamma(1)} \big( {2n+2-n+\frac{n^{2}-(n+1)^{2}}{2}\big)F_{n}
-\big(n+1}+{n^{2}-(n+1)^{2}}\big)F_{n-1}  \]
 \[  u_{n+1}-u_{n}=\\{h} \big( {n+2+\frac{1}{2}n^{2}-\frac{1}{2} n^{2} -n-\frac{1}{2}\big)F_{n}
-\big(n+1}+{n^{2}-n^{2}-n-\frac{1}{2}}\big)F_{n-1}  \]
 \[  u_{n+1}-u_{n}=\\{h} \big( {\frac{3}{2}F_{n}
-\frac{1}{2}}F_{n-1}\big)  \]
\end{remark}
To find the numerical Scheme in the real space we need to apply the inverse Laplace (\ref{eqnum1}). We obtain the following iterative scheme in the real space:  
\begin{eqnarray*}
u_{n+1}(x,t)=\mathcal{L}^{-1}\bigg\{u_n +\frac{h^{\alpha}}{\Gamma(\alpha)} \Big[\bigg( \frac{2 ({n+1})^{\alpha}-n^{\alpha}}{\alpha}+\frac{n^{\alpha+1}-(n+1)^{\alpha+1}}{\alpha+1}\bigg)F_{n}\qquad\\ 
\\-\bigg(\frac{(n+1)^{\alpha}}{\alpha} +\frac{n^{\alpha+1}-(n+1)^{\alpha+1}}{\alpha+1}\bigg)F_{n-1}\Big] \bigg\} 
\end{eqnarray*}
The above equation can be discretised in $x$ using any classical method, including but not limited to forward, backward difference, Crank Nicolson. 
\section{Error Analysis of the Laplace Adam-Bashforth Method}
Let 
\begin{equation}
 {}^{C}_{0}\!\textit{D}^{\alpha} u(x,t) =Lu(x,t)+N u(x,t) 
 \end{equation}
 be a general fractional partial differential equation. As we established earlier the numerical solutions using Laplace Adam-Bashforth method is given as 
 \begin{eqnarray}
u_{n+1}(x,t)=\mathcal{L}^{-1}\bigg\{u_n +\frac{h^{\alpha}}{\Gamma(\alpha)} \Big[\bigg( \frac{2 ({n+1})^{\alpha}-n^{\alpha}}{\alpha}+\frac{n^{\alpha+1}-(n+1)^{\alpha+1}}{\alpha+1}\bigg)F_{n}\qquad\\ \nonumber
\\-\bigg(\frac{(n+1)^{\alpha}}{\alpha} +\frac{n^{\alpha+1}-(n+1)^{\alpha+1}}{\alpha+1}\bigg)F_{n-1}\Big] + R_{n}^{\alpha} \bigg\}\qquad \nonumber
\end{eqnarray} 
Where \[ R_{n}^{\alpha} < \infty  \]

\begin{proof}
Following the derivation presented earlier
 \begin{eqnarray*}   
 u_{n+1}-u_{n}=\frac{1}{\Gamma(\alpha)} \Big[\int_{0}^{t_{n+1}}{(t_{n+1}-\tau)^{\alpha-1}F(u,\tau)d\tau-\int_{0}^{t_{n}} {(t_{n}-\tau)^{\alpha-1}F(u,\tau)d\tau}}\Big]  \end{eqnarray*} 
Where
\[F(u,\zeta)=  \frac{t-t_{n-1}}{t_{n}-t_{n-1}} F_{n}+\frac{t-t_{n}}{t_{n-1}-t_{n}} F_{n-1}+\frac{F^{(2)}(u, \zeta)}{2!}\prod_{i=0}^{1}{(t-t_{i})} \]
\[u_{n+1}-u_{n}=\frac{1}{\Gamma(\alpha)}{\int_{0}^{t_{n+1}}(t_{n+1}-t)^{\alpha-1}\Big(\frac{t-t_{n-1}}{t_{n}-t_{n-1}} F_{n}+\frac{t-t_{n}}{t_{n-1}-t_{n}} F_{n-1}\Big)dt}\]
\[-\frac{1}{\Gamma(\alpha)}{\int_{0}^{t_{n}}(t_{n}-t)^{\alpha-1}\Big(\frac{t-t_{n-1}}{t_{n}-t_{n-1}} F_{n}+\frac{t-t_{n}}{t_{n-1}-t_{n}} F_{n-1}\Big)dt}\]
\[  +\frac{1}{\Gamma(\alpha)}{\int_{0}^{t_{n+1}}\frac{F^{(2)}(u, \zeta)}{2!}\prod_{i=0}^{1}{(t-t_{i})}(t_{n+1}-t)^{\alpha-1}dt}  \]
\[-\frac{1}{\Gamma(\alpha)}{\int_{0}^{t_{n}}\frac{F^{(2)}(u, \zeta)}{2!}\prod_{i=0}^{1}{(t-t_{i})}(t_{n}-t)^{\alpha-1}dt}\]
This is equal to
\begin{eqnarray*} u_{n+1}-u_{n}=\frac{h^{\alpha}}{\Gamma(\alpha)} \Big[\bigg( \frac{2 ({n+1})^{\alpha}-n^{\alpha}}{\alpha}+\frac{n^{\alpha+1}-(n+1)^{\alpha+1}}{\alpha+1}\bigg)F_{n}&\\ 
-\bigg(\frac{(n+1)^{\alpha}}{\alpha} +\frac{n^{\alpha+1}-(n+1)^{\alpha+1}}{\alpha+1}\bigg)F_{n-1} \Big] +R_{n}^{\alpha}  \end{eqnarray*}
Therefore one can easily deduce
\begin{eqnarray*} R_{n}^{\alpha}=\frac{1}{\Gamma(\alpha)}\bigg({\int_{0}^{t_{n+1}}\frac{F^{(2)}(u, \zeta)}{2!}\prod_{i=0}^{1}{(t-t_{i})}(t_{n+1}-t)^{\alpha-1}dt}\\
-{\int_{0}^{t_{n}}\frac{F^{(2)}(u, \zeta)}{2!}\prod_{i=0}^{1}{(t-t_{i})}(t_{n}-t)^{\alpha-1}dt}\bigg)\end{eqnarray*}
\begin{eqnarray*}
|R_{n}^{\alpha}|\le \frac{1}{\Gamma(\alpha)}\bigg({\int_{0}^{t_{n+1}}\bigg|\frac{F^{(2)}(u, \zeta)}{2!}\prod_{i=0}^{1}{(t-t_{i})}(t_{n+1}-t)^{\alpha-1}\bigg|dt}\\
+{\int_{0}^{t_{n}}\bigg|\frac{F^{(2)}(u, \zeta)}{2!}\prod_{i=0}^{1}{(t-t_{i})}(t_{n}-t)^{\alpha-1}dt}\bigg|  \bigg)
\end{eqnarray*}
\begin{eqnarray*}
|R_{n}^{\alpha}|\le \frac{h^{2}}{8\Gamma(\alpha)} \max_{\zeta \in (0, t_{n+1})}\!\!\{{F^{(2)}(u, \zeta)}\}\bigg({\int_{0}^{t_{n+1}}\big|(t_{n+1}-t)^{\alpha-1}\big|dt}\\
+{\int_{0}^{t_{n}}\big|(t_{n}-t)^{\alpha-1}dt}\big|  \bigg)
\end{eqnarray*}
\begin{eqnarray*}
|R_{n}^{\alpha}|\le \frac{h^{2}}{8\Gamma(\alpha)} \max_{\zeta \in (0, t_{n+1})}\!\!\{{F^{(2)}(u, \zeta)}\}\bigg(\frac{t_{n+1}^{\alpha} +t_{n}^{\alpha}}{\alpha} \bigg)\\
\end{eqnarray*}
\begin{eqnarray*}
|R_{n}^{\alpha}|\le \frac{h^{2} h^{\alpha}}{8\Gamma(\alpha+1)} \max_{\zeta \in (0, t_{n+1})}\!\!\{{F^{(2)}(u, \zeta)}\}({n+1}^{\alpha} +{n}^{\alpha}) < + \infty
\end{eqnarray*}
\end{proof}

\section{Applications}

In this section we apply the method to some integer and non integer  order P.D.E. We start with the classical wave equation.
\subsection{Example 1: Integer Order P.D.E}
Consider the following wave-equation like:
\begin{equation}
\frac{\partial u(x,t)}{\partial t}=c ~\frac{\partial u(x,t)}{\partial x}
\end{equation} 
Applying Laplace transform on both sides we get the following
\begin{equation*}
\frac{d u(x,t)}{d t}=c ~(p~u(p,t)-u(0,t))
\end{equation*} 
Silencing the variable $p$ writing $u(p,t)=u(t)$ the later equation can be rewritten as 
\begin{equation*}
\frac{d u(t)}{d t}=F(t, u(t))
\end{equation*} 
The later equation is (\ref{eq2}).\\
We proved earlier that the solution in the Laplace space tho the previous equation is given by (\ref{eq23}).
\begin{equation*}u_{n+1}=u_{n}+h\bigg(\frac{3}{2}F_{n}-\frac{1}{2}F_{n-1}\bigg)  \end{equation*}
Applying the inverse transform to (\ref{eq23}) we have :
\begin{equation}
u(x, t_{n+1})=u(x,t_{n})+\frac{3h}{2}c\frac{\partial u(x, t_{n})}{\partial x}-\frac{h}{2}c\frac{\partial u(x, t_{n-1})}{\partial x}
\end{equation}
Discretising in the space variable we have 
\begin{equation*}
u(x_i, t_{n+1})=u(x_i,t_{n})+\frac{3h}{2}c\bigg[\frac{u(x_{i+1}, t_{n})-u(x_{i}, t_{n})}{\Delta x}\bigg]-\frac{h}{2}c\bigg[\frac{ u(x_{i+1}, t_{n-1})-u(x_{i}, t_{n-1})}{\Delta x}\bigg]
\end{equation*}
Leting $u(x_i, t_{n})=u_i^{n}, \textit{and}~ \Delta x=l$
\begin{equation*}
u_{i}^{n+1}=u_{i}^{n}+\frac{3 }{2} h c \frac{u_{i+1}^{n}-u_{i}^{n}}{l}-\frac{1}{2}hc \frac{ u_{i+1}^{n-1}-u_{i}^{n-1}}{l}
\end{equation*}
\begin{equation}
u_{i}^{n+1}=(1-\frac{3 h c}{2l}  )u_{i}^{n}+\frac{3 h c}{2l}u_{i+1}^{n} -\frac{hc}{2l} u_{i+1}^{n-1}+\frac{hc}{2l}u_{i}^{n-1} \label{eq53}
\end{equation}
\subsubsection{Stability Analysis}

Assume we have a Fourier expansion in space of \[u(x,t)=\sum_{f}\widehat{u}(t)\exp(jfx) \]
Equation (\ref{eq53}) becomes 
\begin{equation*}
\widehat{u}_{n+1}e^{jifl}=(1-\frac{3 h c}{2l} )\widehat{u}_{n} e^{jifl}+\frac{3 h c}{2l}\widehat{u}_{n}e^{j(i+1)fl} -\frac{hc}{2l} \widehat{u}_{n-1}e^{j(i+1)fl}+\frac{hc}{2l}\widehat{u}_{n-1}e^{jifl}
\end{equation*}
this is 
\begin{equation*}
\widehat{u}_{n+1}=(1-\frac{3 h c}{2l} )\widehat{u}_{n} +\frac{3 h c}{2l}\widehat{u}_{n}e^{jfl} -\frac{hc}{2l} \widehat{u}_{n-1}e^{jfl}+\frac{hc}{2l}\widehat{u}_{n-1}
\end{equation*}

\begin{equation*}
\widehat{u}_{n+1}=(1-\frac{3 h c}{2l}  +\frac{3 h c}{2l}e^{jfl})\widehat{u}_{n} +(\frac{hc}{2l}-\frac{hc}{2l} e^{jfl})\widehat{u}_{n-1}
\end{equation*}
\begin{equation*}
\frac{\widehat{u}_{n+1}}{\widehat{u}_{n}}=(1-\frac{3 h c}{2l} ) +\frac{3 h c}{2l}e^{jfl} -\frac{hc}{2l} \frac{\widehat{u}_{n-1}}{\widehat{u}_{n}}e^{jfl}+\frac{hc}{2l}\frac{\widehat{u}_{n-1}}{\widehat{u}_{n}}
\end{equation*}

Finding the ration we have:
\begin{equation*}
\frac{\widehat{u}_{n+1}}{\widehat{u}_{n} }=1-\frac{3 h c}{2l}  +\frac{3 h c}{2l}e^{jfl}+(\frac{hc}{2l}-\frac{hc}{2l} e^{jfl})\frac{\widehat{u}_{n-1}}{\widehat{u}_{n} }
\end{equation*}
\begin{equation*}
\frac{\widehat{u}_{n+1}}{\widehat{u}_{n} }=1-\frac{3 h c}{2l}\big(1-cos(fl)\big) +j\frac{3 h c}{2l}sin{(fl)} +\frac{hc}{2l}(1-\cos{(jfl)})\frac{\widehat{u}_{n-1}}{\widehat{u}_{n} }-j\frac{hc}{2l}\sin{(fl)}\frac{\widehat{u}_{n-1}}{\widehat{u}_{n} }
\end{equation*}
\begin{equation*}
\frac{\widehat{u}_{n+1}}{\widehat{u}_{n} }=1-\frac{3 h c}{l}\sin^2\bigg(\frac{fl}{2}\bigg) +\frac{hc}{l}\sin^2\bigg(\frac{fl}{2}\bigg)\frac{\widehat{u}_{n-1}}{\widehat{u}_{n} }+j\frac{h c}{2l}sin{(fl)}(3 -\frac{\widehat{u}_{n-1}}{\widehat{u}_{n} })
\end{equation*}

Let us prove that $\forall ~n~ |u_n|< |u_0|.$
For $n=0$ 
\[ \bigg|\frac{\widehat{u}_1}{\widehat{u}_0}\bigg|=\bigg| 1-\frac{3 h c}{2l}\sin^2\bigg(\frac{fl}{2}\bigg)\bigg|   \]

This means that we have 

\[ \bigg|\frac{\widehat{u}_1}{\widehat{u}_0}\bigg|<1 \Longleftrightarrow -1< 1-\frac{3 h c}{2l}\sin^2\bigg(\frac{fl}{2}\bigg)<1   \]
\[ \Longrightarrow -2< -\frac{3 h c}{2l}\sin^2\bigg(\frac{fl}{2}\bigg)<0   \]
\[ \Longrightarrow 0< \frac{3 h c}{2l}\sin^2\bigg(\frac{fl}{2}\bigg)<2   \]
\[ \Longrightarrow 0< \frac{3 h c}{4l}\sin^2\bigg(\frac{fl}{2}\bigg)<1   \]
The previous condition will certainly be achieved if 
\[ 0< \frac{3 h c}{4l}<1   \] since

\[  0< \frac{3 h c}{4l}\sin^2\bigg(\frac{fl}{2}\bigg)\le \frac{3 h c}{4l} .  \]
\[\frac{3 h c}{4l}<1  \Longrightarrow \]
\begin{equation} 
\frac{h}{l}<\frac{4}{3}c  \label{cond} 
\end{equation}.\\\\

Now let us assume $ \forall j \le n, |\widehat{u}_n|<|\widehat{u}_0|$ and prove that $\widehat{u}_{n+1}<|\widehat{u}_0|$
The condition (\ref{cond}) ensures that the coefficients in the following equation obtained earlier above are positive.
\[\widehat{u}_{n+1}=(1-\frac{3 h c}{2l} )\widehat{u}_{n} +\frac{3 h c}{2l}\widehat{u}_{n}e^{jfl} -\frac{hc}{2l} \widehat{u}_{n-1}e^{jfl}+\frac{hc}{2l}\widehat{u}_{n-1}\].
The later then implies that 
\[|\widehat{u}_{n+1}|<(1-\frac{3 h c}{2l} )|\widehat{u}_{n}| +\frac{3 h c}{2l}|\widehat{u}_{n}||e^{jfl}| -\frac{hc}{2l} |\widehat{u}_{n-1}||e^{jfl}|+\frac{hc}{2l}|\widehat{u}_{n-1}|\].
Using the induction hypothesis this becomes 
\[|\widehat{u}_{n+1}|<(1-\frac{3 h c}{2l} )|\widehat{u}_{0}| +\frac{3 h c}{2l}|\widehat{u}_{0}| -\frac{hc}{2l} |\widehat{u}_{0}|+\frac{hc}{2l}|\widehat{u}_{0}|\].
and therefore we obtain
\[|\widehat{u}_{n+1}|<\bigg(1-\frac{3 h c}{2l}  +\frac{3 h c}{2l} -\frac{hc}{2l} +\frac{hc}{2l}\bigg)|\widehat{u}_{0}|\].
\[|\widehat{u}_{n+1}|<|\widehat{u}_{0}|\].
We can therefore conclude that the numerical scheme solution presented above is stable.
\subsection{Example: Fractional Order P.D.E}
Consider the following fractional Order P.D.E, the fractional derivative is given in the Caputo Sense.
\begin{equation}
\frac{\partial^{\alpha} u(x,t)}{\partial t^{\alpha}}=d~\frac{\partial^2 u(x,t)}{\partial x^2}
\end{equation} 
Applying Laplace transform on both sides we get the following
\begin{equation*}
{}_{0}^{C}\!\!\textit{D}_{t}^{\alpha} {u(p,t)}=d~ \mathcal{L}\Big\{\frac{\partial^2 u(x,t)}{\partial x^2}\Big\}
\end{equation*}
 
\begin{equation*}
~~~~~~~~~~~~~~~~~~~~~~~~=d~ \mathcal{L}\Big(p^2 u(p,t)-p~u(0,t) -u(0,t)\Big)
\end{equation*} 
Silencing the variable $p$ writing $u(p,t)=u(t)$ the later equation can be rewritten as 
\begin{equation*}
{}_{0}^{C}\!\!\textit{D}_{t}^{\alpha} {u(p,t)}=d~ F(t, u(t))
\end{equation*}
The later equation is (\ref{eq32}).\\
We proved earlier that the solution in the Laplace space tho the previous equation is given by (\ref{eqnum1})
\begin{eqnarray*} u_{n+1}-u_{n}=\frac{h^{\alpha}}{\Gamma(\alpha)} \Big[\bigg( \frac{2 ({n+1})^{\alpha}-n^{\alpha}}{\alpha}+\frac{n^{\alpha+1}-(n+1)^{\alpha+1}}{\alpha+1}\bigg)F_{n}&\\ 
-\bigg(\frac{(n+1)^{\alpha}}{\alpha} +\frac{n^{\alpha+1}-(n+1)^{\alpha+1}}{\alpha+1}\bigg)F_{n-1} \Big]
\end{eqnarray*} 
Applying the inverse transform to (\ref{eqnum1}) we have :
\begin{equation*}
u(x, t_{n+1})=u(x,t_{n})+\frac{h^{\alpha}}{\Gamma(\alpha)} \delta_{n}^{\alpha} d \frac{\partial^2 u(x, t_{n})}{\partial x^2}-\frac{h^{\alpha}}{\Gamma{(\alpha)}} \delta_{n}^{\alpha,1} d \frac{\partial^2 u(x, t_{n-1})}{\partial x^2}
\end{equation*}
where 
\[ \delta_{n}^{\alpha} = \frac{2 ({n+1})^{\alpha}-n^{\alpha}}{\alpha}+\frac{n^{\alpha+1}-(n+1)^{\alpha+1}}{\alpha+1} \]and 
\[\delta_{n}^{\alpha,1}=\frac{(n+1)^{\alpha}}{\alpha} +\frac{n^{\alpha+1}-(n+1)^{\alpha+1}}{\alpha+1}    \]
Discretising in the space variable we have 
\begin{eqnarray*}
u(x_i, t_{n+1})=u(x_i,t_{n})+\frac{h^{\alpha}}{\Gamma(\alpha)} \delta_{n}^{\alpha} d \frac{ u(x_{i+1}, t_{n})-2 u(x_i, t_n)+u(x_{i-1}, t_{n})}{(\Delta x)^2}&\\
-\frac{h^{\alpha}}{\Gamma{(\alpha)}} \delta_{n}^{\alpha,1} d \frac{u(x_{i+1}, t_{n-1})-2 u(x_i, t_{n-1})+u(x_{i-1}, t_{n-1})}{(\Delta x)^2}
\end{eqnarray*}
Leting $u(x_i, t_{n})=u_i^{n}, \textit{and}~ \Delta x=l$ the equation becomes 
\begin{eqnarray*}
u_i^{n+1}=u_i^{n}+\frac{h^{\alpha}}{\Gamma(\alpha)} \delta_{n}^{\alpha} d \frac{ u_{i+1}^{n}-2 u_i^{n}+u_{i-1}^{n}}{l^2}&\\
-\frac{h^{\alpha}}{\Gamma{(\alpha)}} \delta_{n}^{\alpha,1} d \frac{u_{i+1}^{n-1}-2 u_i^{n-1}+u_{i-1}^{n-1}}{l^2}
\end{eqnarray*}
This is
\begin{eqnarray}
u_i^{n+1}=\Big(1-2\frac{h^{\alpha} d}{l^2\Gamma(\alpha)} \delta_{n}^{\alpha}  \Big)u_i^{n}+  \frac{h^{\alpha} d}{l^2\Gamma(\alpha)} \delta_{n}^{\alpha} u_{i+1}^{n}+\frac{h^{\alpha} d}{l^2\Gamma(\alpha)} \delta_{n}^{\alpha}u_{i-1}^{n}\nonumber &\\ 
-\frac{h^{\alpha} d}{l^2\Gamma{(\alpha)}} \delta_{n}^{\alpha,1} u_{i+1}^{n-1}+2\frac{h^{\alpha} d}{l^2\Gamma{(\alpha)}} \delta_{n}^{\alpha,1} u_i^{n-1}-\frac{h^{\alpha} d}{l^2\Gamma{(\alpha)}} \delta_{n}^{\alpha,1}u_{i-1}^{n-1} \label{eq55}
\end{eqnarray}
\subsubsection{Stability Analysis of the numerical scheme solution fractional Order P.D.E}
Assume we have a Fourier expansion in space of \[u(x,t)=\sum_{f}\widehat{u}(t)\exp(jfx) \] while leting \[u_i^ n=\widehat{u}_n\exp{(jfi\Delta x)}=\widehat{u}_n\exp{(jfil)}   \]
Equation (\ref{eq55}) becomes
\begin{eqnarray*}
\widehat{u}_{n+1} e^{jfil}=\Big(1-2\frac{h^{\alpha} d}{l^2\Gamma(\alpha)} \delta_{n}^{\alpha}  \Big)\widehat{u}_{n}e^{jfil}+  \frac{h^{\alpha} d}{l^2\Gamma(\alpha)} \delta_{n}^{\alpha} \widehat{u}_{n}e^{jf(i+1)l}+\frac{h^{\alpha} d}{l^2\Gamma(\alpha)} \delta_{n}^{\alpha}\widehat{u}_{n} e^{jf(i-1)l}&\\
-\frac{h^{\alpha} d}{l^2\Gamma{(\alpha)}} \delta_{n}^{\alpha,1} \widehat{u}_{n-1}e^{jf(i+1)l}+2\frac{h^{\alpha} d}{l^2\Gamma{(\alpha)}} \delta_{n}^{\alpha,1} \widehat{u}_{n-1}e^{jfil}-\frac{h^{\alpha} d}{l^2\Gamma{(\alpha)}} \delta_{n}^{\alpha,1}\widehat{u}_{n-1}e^{jf(i-1)l} 
\end{eqnarray*} 
This is 
 \begin{eqnarray*}
\widehat{u}_{n+1} =\Big(1-2\frac{h^{\alpha} d}{l^2\Gamma(\alpha)} \delta_{n}^{\alpha}  \Big)\widehat{u}_{n}+  \frac{h^{\alpha} d}{l^2\Gamma(\alpha)} \delta_{n}^{\alpha} \widehat{u}_{n}e^{jfl}+\frac{h^{\alpha} d}{l^2\Gamma(\alpha)} \delta_{n}^{\alpha}\widehat{u}_{n} e^{-jfl}&\\
-\frac{h^{\alpha} d}{l^2\Gamma{(\alpha)}} \delta_{n}^{\alpha,1} \widehat{u}_{n-1}e^{jfl}+2\frac{h^{\alpha} d}{l^2\Gamma{(\alpha)}} \delta_{n}^{\alpha,1} \widehat{u}_{n-1}-\frac{h^{\alpha} d}{l^2\Gamma{(\alpha)}} \delta_{n}^{\alpha,1}\widehat{u}_{n-1}e^{-jfl} 
\end{eqnarray*}
Which implies 
  \begin{eqnarray*}
\widehat{u}_{n+1} =\Big(1-2\frac{h^{\alpha} d}{l^2\Gamma(\alpha)} \delta_{n}^{\alpha}  \Big)\widehat{u}_{n}+2\frac{h^{\alpha} d}{l^2\Gamma{(\alpha)}} \delta_{n}^{\alpha,1} \widehat{u}_{n-1}+  \frac{h^{\alpha} d}{l^2\Gamma(\alpha)} \delta_{n}^{\alpha} \widehat{u}_{n}\big(e^{jfl}+ e^{-jfl}\big)&\\
-\frac{h^{\alpha} d}{l^2\Gamma{(\alpha)}} \delta_{n}^{\alpha,1} \widehat{u}_{n-1}\big(e^{jfl}+e^{-jfl}\big) 
\end{eqnarray*}
Further we have 
  \begin{eqnarray*}
\widehat{u}_{n+1} =\Big(1-2\frac{h^{\alpha} d}{l^2\Gamma(\alpha)} \delta_{n}^{\alpha}  \Big)\widehat{u}_{n}+2\frac{h^{\alpha} d}{l^2\Gamma{(\alpha)}} \delta_{n}^{\alpha,1} \widehat{u}_{n-1}+  \frac{h^{\alpha} d}{l^2\Gamma(\alpha)} \delta_{n}^{\alpha} \widehat{u}_{n}\big(2\cos{(fl)}\big)&\\
-\frac{h^{\alpha} d}{l^2\Gamma{(\alpha)}} \delta_{n}^{\alpha,1} \widehat{u}_{n-1}\big(2\cos{(fl)}\big) 
\end{eqnarray*}
and 
 \begin{eqnarray*}
\widehat{u}_{n+1} =\Big(1-\frac{2 h^{\alpha} d}{l^2\Gamma(\alpha)} \delta_{n}^{\alpha}+\frac{2 h^{\alpha} d}{l^2\Gamma(\alpha)} \delta_{n}^{\alpha}\cos{(fl)}  \Big)\widehat{u}_{n}+\Big(\frac{2 h^{\alpha} d}{l^2\Gamma{(\alpha)}} \delta_{n}^{\alpha,1}&\\
-\frac{2 h^{\alpha} d}{l^2\Gamma{(\alpha)}} \delta_{n}^{\alpha,1} \cos{(fl)}  \Big)\widehat{u}_{n-1} 
\end{eqnarray*}
Thus we can write
 \begin{eqnarray}
\widehat{u}_{n+1} =\Big(1-\frac{2 h^{\alpha} d}{l^2\Gamma(\alpha)} \delta_{n}^{\alpha}\big(1-\cos{(fl)}\big)  \Big)\widehat{u}_{n}+\Big(\frac{2 h^{\alpha} d}{l^2\Gamma{(\alpha)}} \delta_{n}^{\alpha,1}\big(1-cos{(fl)}\big)  \Big)\widehat{u}_{n-1} ~~\label{eq56}
\end{eqnarray}
From equation (\ref{eq56}) we have 
 \begin{eqnarray*}
\frac{\widehat{u}_{n+1}}{\widehat{u}_{n}} =1-\frac{2 h^{\alpha} d}{l^2\Gamma(\alpha)} \delta_{n}^{\alpha}\big(1-\cos{(fl)}\big)+\Big(\frac{2 h^{\alpha} d}{l^2\Gamma{(\alpha)}} \delta_{n}^{\alpha,1}\big(1-cos{(fl)}\big)\frac{\widehat{u}_{n-1}}{\widehat{u}_n} \Big)
\end{eqnarray*}
 \begin{eqnarray*}
~=1-\frac{4 h^{\alpha} d}{l^2\Gamma(\alpha)} \delta_{n}^{\alpha}\sin^2{\bigg(\frac{fl}{2}\bigg)} +\Big(\frac{4 h^{\alpha} d}{l^2\Gamma{(\alpha)}} \delta_{n}^{\alpha,1} \sin^2\!\!{\bigg(\frac{fl}{2}\bigg)} \frac{\widehat{u}_{n-1}}{\widehat{u}_n} \Big)
\end{eqnarray*}
For $n=0$
\[\frac{\widehat{u}_1}{\widehat{u}_0}= 1-\frac{4 h^{\alpha} d}{l^2\Gamma(\alpha)} \delta_{n}^{\alpha}\sin^2\!\!{\bigg(\frac{fl}{2}\bigg)}  \]

\[\bigg|\frac{\widehat{u}_1}{\widehat{u}_0}\bigg|<1 \Longleftrightarrow -1<1-\frac{4 h^{\alpha} d}{l^2\Gamma(\alpha)} \delta_{n}^{\alpha}\sin^2\!\!{\bigg(\frac{fl}{2}\bigg)} <1   \Longrightarrow    \]

\[0<\frac{4 h^{\alpha} d}{l^2\Gamma(\alpha)} \delta_{n}^{\alpha}\sin^2\!\!{\bigg(\frac{fl}{2}\bigg)} <2   \Longrightarrow 0<\frac{2 h^{\alpha} d}{l^2\Gamma(\alpha)} \delta_{n}^{\alpha}\sin^2\!\!{\bigg(\frac{fl}{2}\bigg)} <1    \]

This will certainly be achieved if 

\[0<\frac{2 h^{\alpha} d}{l^2\Gamma(\alpha)} \delta_{n}^{\alpha} <1    \]
as 
\[\frac{2 h^{\alpha} d}{l^2\Gamma(\alpha)} \delta_{n}^{\alpha}\sin^2\!\!{\bigg(\frac{fl}{2}\bigg)}\le\frac{2 h^{\alpha} d}{l^2\Gamma(\alpha)} \delta_{n}^{\alpha}.   \]

Therefore  \[0<\frac{2 h^{\alpha} d}{l^2\Gamma(\alpha)} \delta_{n}^{\alpha} <1   \Longrightarrow \]
\begin{equation}
\frac{ h^{\alpha} }{l^2} <\frac{\Gamma(\alpha)}{2d \delta_{n}^{\alpha}}. \label{eq57}
\end{equation}

We will now prove that $\forall~ n~ |u_n| <|u_0|$, we proved already that $|u_1|<|u_0|$, let us assume that $|u_j|<|u_0|, \forall j\le n $ and prove that $|u_{n+1}|<|u_0|$
   \begin{eqnarray*}
|\widehat{u}_{n+1}|=\bigg|\bigg(1-\frac{4 h^{\alpha} d}{l^2\Gamma(\alpha)} \delta_{n}^{\alpha}\sin^2{\big(\frac{fl}{2}\big)\bigg)}\widehat{u}_n +\Big(\frac{4 h^{\alpha} d}{l^2\Gamma{(\alpha)}} \delta_{n}^{\alpha,1} \sin^2\!\!{\bigg(\frac{fl}{2}\bigg)}\Big)\widehat{u}_{n-1}\bigg|
\end{eqnarray*}
\[Let ~ A_1=\bigg(1-\frac{4 h^{\alpha} d}{l^2\Gamma(\alpha)} \delta_{n}^{\alpha}\sin^2{\big(\frac{fl}{2}\big)\bigg)} , A_2=\Big(\frac{4 h^{\alpha} d}{l^2\Gamma{(\alpha)}} \delta_{n}^{\alpha,1} \sin^2\!\!{\bigg(\frac{fl}{2}\bigg)}\Big)   \]
Because of the condition (\ref{eq57}) we have both  $ A_1, A_2>0.$ 
\[ |\widehat{u}_{n+1}|=\bigg|\bigg(1-\frac{4 h^{\alpha} d}{l^2\Gamma(\alpha)} \delta_{n}^{\alpha}\sin^2{\big(\frac{fl}{2}\big)\bigg)}\widehat{u}_n +\Big(\frac{4 h^{\alpha} d}{l^2\Gamma{(\alpha)}} \delta_{n}^{\alpha,1} \sin^2\!\!{\bigg(\frac{fl}{2}\bigg)}\Big)\widehat{u}_{n-1}\bigg|\]
\[=A_1|\widehat{u}_n|+A_2|\widehat{u}_{n-1}|\]
\[ |\widehat{u}_{n+1}|< A_1|\widehat{u}_n|+A_2|\widehat{u}_{n-1}|\]
By induction hypothesis we have 
\[ |\widehat{u}_{n+1}|<|A_1||\widehat{u}_0|+|A_2||\widehat{u}_0| \]
This means
\[ |\widehat{u}_{n+1}|<(A_1+A_2)|\widehat{u}_0| \]
But 
\[A_1+A_2=1-\frac{4 h^{\alpha} d}{l^2\Gamma(\alpha)} \delta_{n}^{\alpha}\sin^2{\big(\frac{fl}{2}\big)}+ \frac{4 h^{\alpha} d}{l^2\Gamma{(\alpha)}} \delta_{n}^{\alpha,1} \sin^2\!\!{\bigg(\frac{fl}{2}\bigg)} =1\]
Therefore
\[ |\widehat{u}_{n+1}|<|\widehat{u}_0| \]
This prove that $\forall~ n~ |u_n| <|u_0|$ the numerical scheme solution to the  fractional Order P.D.E is stable.\\

\section{Graphical Simulations}
\newpage
\begin{figure}
	\centering\includegraphics[width=0.9\linewidth]{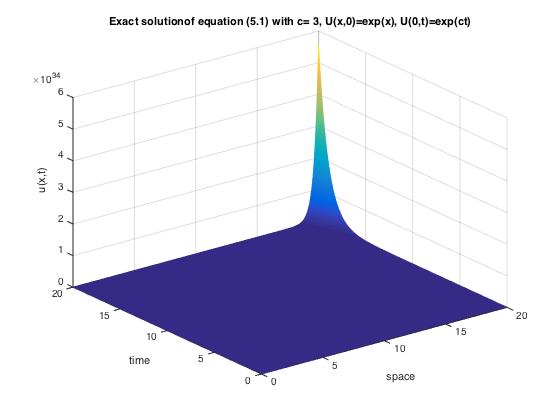}
	\caption{Exact Solution of equation (5.1) with u(x,0)=exp(x); u(0,t)=exp(ct).}
\end{figure}
\begin{figure}
	\centering\includegraphics[width=0.9\linewidth]{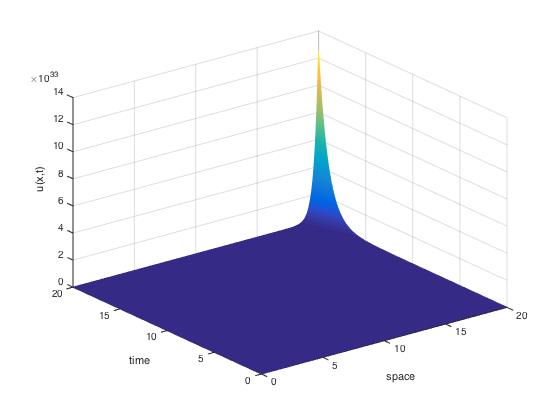}
	\caption{Approximate Solution of equation (5.1) given by (5.3), with u(x,0)=exp(x); u(0,t)=exp(ct).}
\end{figure}
\begin{figure}
	\centering\includegraphics[width=0.9\linewidth]{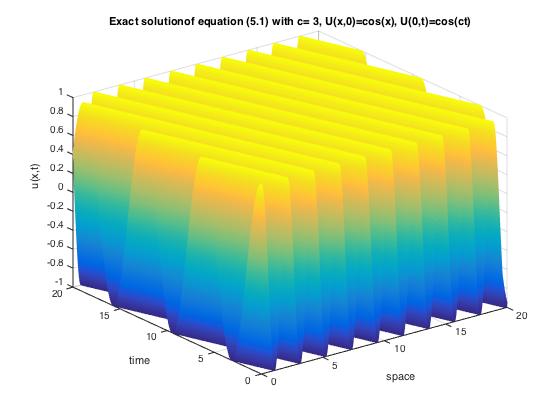}
	\caption{Exact Solution of equation (5.1) with u(x,0)=cos(x); u(0,t)=cos(ct).}
\end{figure}
\begin{figure}
	\centering\includegraphics[width=0.9\linewidth]{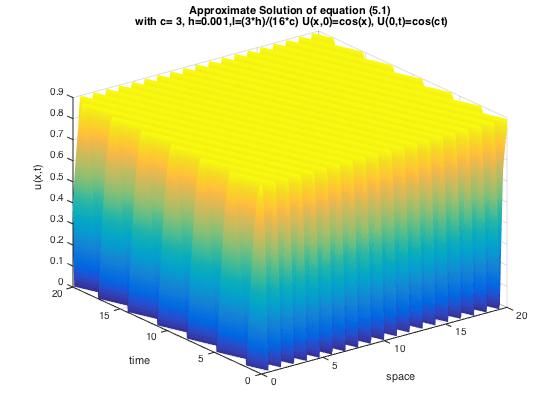}
	\caption{Approximate Solution of equation (5.1) given by (5.3), with u(x,0)=cos(x); u(0,t)=cos(ct).}
\end{figure}

\newpage
\section{Conclusion}
In order to extend the well-known Adams-Bashforth numerical scheme to partial differential equations with integer and non-integer order derivatives, we introduced a new reliable and efficient numerical scheme. The method is a combination of Laplace transform, Adams-Bashforth and (forward or backward numerical scheme).  We developed the method for general partial differential equations with local and non-local differentiation. We presented in detail the error analysis and the convergence of the method. In the case of fractional partial differential equations, the method provides a numerical algorithm that is easier to implement. Unlike the conventional methods, forward, backward, Crank-Nicholson, the cumbersome summation that always appears in the additional term of their numerical algorithm for the case of fractional partial differential equations, does not exist with our method. This leads to an easier proof of stability and convergence. We illustrated the method by solving two partial differential equations including wave equation for the local case and a diffusion equation for fractional case. We studied the stability of each example. The proof shows without doubt that our method is very stable and also converges very quickly to the exact solution. We believe this method will turn out to be a very useful numerical scheme that will help solving nonlinear and linear partial differential equations with local and non-local operators.
\newpage
\section*{Conflict of Interests Disclosure}
The authors declare no conflict of Interests.

\end{document}